\documentclass[a4paper,11pt]{amsart}
\pdfoutput=1

\usepackage{a4wide}
\usepackage{amsfonts}
\usepackage{amscd}
\usepackage{amssymb}
\usepackage{amsthm}
\usepackage{amsmath}
\usepackage{graphicx}
\usepackage{subfigure}
\usepackage{dsfont}
\usepackage[american]{babel}
\usepackage{algorithmic}
\usepackage[geometry]{ifsym}

\newtheorem{theorem}{Theorem}
\newtheorem{lemma}[theorem]{Lemma}

\theoremstyle{definition}

\newtheorem{algorithm}[theorem]{Algorithm}

\theoremstyle{remark}
\newtheorem{remark}[theorem]{Remark}

\newcommand{\N}{{\mathbb N}}
\newcommand{\R}{{\mathbb R}}

\newcommand{\cF}{{\mathcal F}}

\newcommand{\cK}{{\mathcal K}}

\newcommand{\cL}{{\mathcal L}}

\newcommand{\cN}{{\mathcal N}}

\newcommand{\identity}{\mathds{1}}

\newcommand{\set}[1]{{\left\{#1\right\}}}

\newcommand{\size}[1]{\left|#1\right|}

\newcommand{\expct}{{\mathbb E}}

\newcommand{\sip}[2]{\left( #1,#2 \right)}

\newcommand{\integral}[3]{\int_{#1} #2 \; d#3}
\renewcommand{\Re}{\mathrm{Re}\,}
\renewcommand{\Im}{\mathrm{Im}\,}

\newcommand{\fourier}{{\mathcal F}}
\newcommand{\ifourier}{{\mathcal F}^{-1}}
\newcommand{\fft}{{\rm FFT}}
\newcommand{\ifft}{{\rm FFT}^{-1}}

\newcommand{\ie}{i.e.\ }

\newcommand{\bigtimes}{\mathop{\textifgeo{\BigCross}}
\limits}
\newcommand{\Op}{\operatorname{O}}

\title{Fast simulation of Gaussian random fields}

\author[A.~Lang]{Annika Lang}
\address{Annika Lang\newline
Fakult\"at f\"ur Mathematik und Informatik, Universit\"at Mannheim\newline
D--68131 Mann\-heim, Germany\newline
and\newline
Seminar f\"ur Angewandte Mathematik\newline
ETH Z\"urich\newline
CH--8092 Z\"urich, Switzerland}
\email{annika.lang@sam.math.ethz.ch}

\author[J.~Potthoff]{J\"urgen Potthoff}
\address{J\"urgen Potthoff\newline
Fakult\"at f\"ur Mathematik und Informatik, Universit\"at Mannheim\newline
D--68131 Mann\-heim, Germany}
\email{potthoff@math.uni-mannheim.de}

\date{February 11, 2013}

\begin{document}

\begin{abstract}
Fast Fourier transforms are used to develop algorithms for the fast generation of
correlated Gaussian random fields on rectangular regions of $\R^d$. The complexities
of the algorithms are derived, simulation results and error analysis are presented.
\end{abstract}

\maketitle

\section{Introduction}
In this note we present two algorithms for the fast simulation of Gaussian random
fields (GRF) on a rectangular region of $\R^d$.

Among the many possible applications we only want to mention here the simulation of numerical
solutions of stochastic partial differential equations (SPDEs), and this was one of the
original motivations of the present note~\cite{Lang:2007}. A typical example of
an SPDE is the stochastic heat equation
\begin{equation*}
	\frac{\partial}{\partial t}\,u(t,x) = \Delta_x u(t,x) + \eta(t,x),
			\qquad t\ge 0,\,x\in D,
\end{equation*}
in a domain $D$ of $\R^d$, driven by a centered Gaussian random field $\eta(t,x)$,
which typically is ``white'', i.e., uncorrelated, in the time direction, and which
has a covariance function $C(x,y)$, $x$, $y\in D$, in the spatial variable. Here
$\Delta_x$ stands for the Laplacian in the $x$ variable. For accounts of the
mathematical framework for SPDEs we refer the interested reader to, e.g.,
\cite{Chow:2007, DaPrato:1992, Prevot:2007} and to the works cited there.

In this paper we describe algorithms which generate approximations to centered
Gaussian random fields on a rectangular grid with the help of the fast Fourier
transform~(FFT). This entails that the resulting generated samples have periodic
boundary conditions. If one uses instead discrete cosine and sine transforms (DCT
and DST), similar results are obtained with Neumann and Dirichlet boundary
conditions. A generalization to other domains than rectangles should be possible via
NFFT (cf.~\cite{Fourmont:1999, Potts:2003, Fenn:2005}). An algorithm using FFT for
fast Gaussian random field generation for the one-dimensional case can be found
in~\cite{Ripley:1987}. Here we generalize this algorithm to $d$-dimensional
rectangular regions. Furthermore we present a second algorithm that replaces one of
the FFTs by a special data arrangement. In computer experiments this second
algorithm is faster by a factor two. On an older computer the generation of one GRF
of size $512\times 512$ takes $0.25$ seconds ($250 000$ clocks).

Spectral methods of different type can be found in~\cite{Kramer:2007} and references
therein. Here we take the spectral representation of a given covariance, generate
the Gaussian random field by multiplying the spectral density with the Fourier
transformed white noise field and then apply the inverse Fourier transform.
For an approach based on the circular embedding technique we refer the interested
reader to~\cite{Gneiting:2006} and the literature quoted there.

The paper is organized in the following way: Section~\ref{section_GRF} presents the
properties of GRFs that will be used for the construction of the algorithms. These
algorithms are presented in Section~\ref{section_algorithm} and their complexity is
analyzed. Finally Section~\ref{section_tests} presents a class of covariance
functions and shows simulation results and error estimates of an implementation in
C++.

\section{Properties of continuous GRFs}\label{section_GRF}

Assume that $\varphi$ is a real-valued stationary Gaussian random field (GRF) on
$\R^d$. Thus for all $x \in \R^d$, $\varphi(x)$ defines a random variable, and the
family $(\varphi(x),x \in \R^d)$ consists of identically distributed random
variables on a probability space $(\Omega, \cF,P)$ that satisfy for any subset
$\set{x_1, \ldots, x_n \in \R^d}$ and $\alpha_k \in \R$ that the random variable
$\sum_{k=1}^n \alpha_k \, \varphi(x_k)$ is normally distributed. Then its mean $m =
\expct(\varphi(x)) = \expct(\varphi(0))$, $x \in \R^d$, and its covariance $C(x) =
\expct(\varphi(0) \varphi(x)) - m^2 = \expct(\varphi(y) \varphi(x+y))- m^2$, $x, y
\in \R^d$ completely characterize the random field and are therefore its only
relevant statistical parameters. Without loss of generality we will assume that $m =
0$. The goal of this section is to construct GRFs $\varphi$ for a given covariance
$C$. We will present a heuristic approach here. A mathematically rigorous approach
can be found in~\cite{Lang:2007}.

Let us start with the properties of the covariance function. We note that $C$ is a
positive semi-definite, symmetric function, and from now on we assume that it is
continuous. Then Bochner's theorem~\cite{Bochner:1955} states that $C$ is the
Fourier transform of a positive measure $\mu_C$ on $\R^d$, \ie $C$ can be written as
\[
C(x,y) = \int e^{-2\pi i \sip{p}{x-y}} \, d\mu_C(p), \qquad x,y \in \R^d,
\]
where $(\cdot,\cdot)$ denotes the Euclidean inner product on $\R^d$. We observe that
$C$ is actually a function of the difference $x-y$, and that it is an even function
of $x-y$. For most practical purposes, there is no loss of generality, if we assume
that $\mu_C$ has a Lebesgue density denoted by $\gamma$ which is even and positive.
Then $C$ can be written in the following way:
\[
C(x,y) = \integral{\R^d}{e^{-2\pi i \sip{p}{x-y}} \gamma(p)}{p}, \qquad x, y \in
\R^d.
\]
Let $W$ be a Gaussian white noise random field on $\R^d$ also known as cylindrical
Wiener process or $Q$--Wiener process with $Q = \identity$, \ie informally, $W$ is a
centered Gaussian family $\set{W(x), x \in \R^d}$ with covariance $\expct(W(x)
W(y))= \delta(x-y)$, $x,y \in \R^d$. Set
\begin{equation}\label{eqn_first_GRF_alg}
\varphi(x) = ( \ifourier \gamma^{1/2} \fourier W)(x), \qquad x \in \R^d,
\end{equation}
where $\fourier$ denotes the $d$-dimensional Fourier transform and $\ifourier$ is
its inverse. Then, since $W$ is centered Gaussian, so is $\varphi$. The covariance
of $\varphi$ is given by
\begin{align*}
\begin{split}
\expct(\varphi(x) & \varphi(y))\\
& = \iiiint e^{-2\pi i (\sip{p}{x} + \sip{q}{y})} \gamma(p)^{1/2} \gamma(q)^{1/2} e^{2\pi i (\sip{p}{x^\prime} + \sip{q}{y^\prime})} \expct(W(x^\prime)W(y^\prime)) \, dx^\prime \, dy^\prime \, dp \, dq\\
& = \iint e^{-2\pi i (\sip{p}{x} + \sip{q}{y})} \gamma(p)^{1/2} \gamma(q)^{1/2} \int e^{2\pi i \sip{p+q}{x^\prime}} \, dx^\prime \, dp \, dq\\
& = \int e^{-2\pi i \sip{p}{x-y}} \gamma(p) \, dp\\
& = C(x,y).
\end{split}
\end{align*}
Equation~\eqref{eqn_first_GRF_alg} suggests our first algorithm to generate samples
of $\varphi$ with given covariance $C$, and it is the first one presented in
Section~\ref{section_algorithm}.

Next we will work out a more efficient algorithm by replacing one Fourier transform
with a faster operation. Our computer experiments with this algorithm took half the
time of those done with the algorithm based on Equation~\eqref{eqn_first_GRF_alg}.
The goal is to construct $\fourier W$ directly, \ie a complex GRF with the same
distribution as $\fourier W$ has to be generated. We consider the complex-valued
random field $\fourier W$ with
\[
\fourier W(p) = \integral{\R^d}{e^{2\pi i \sip{p}{x}} W(x)}{x}
\]
for $p \in \R^d$. This is obviously centered Gaussian and the covariance is given by
\begin{align*}
\expct(\fourier W(p) \overline{\fourier W(q)})
 = \iint e^{2\pi i (\sip{p}{x} - \sip{q}{y})} \expct(W(x) W(y)) \, dx \, dy
 = \int e^{2\pi i \sip{p-q}{x}} \, dx
 = \delta(p-q),
\end{align*}
and similarly
\begin{align*}
\expct(\fourier W(p) \fourier W(q)) = \delta(p+q).
\end{align*}
The following lemma gives a direct construction of $\fourier W$.

\begin{lemma}\label{lemma_FW}
Let $V$ be the complex-valued random field defined by
\begin{align*}
\begin{split}
\Re V  = \pi^+ W \quad \text{and} \quad
\Im V   = \pi^- W,
\end{split}
\end{align*}
where
\begin{align*}
\begin{split}
\pi^+ W(p) = \tfrac{1}{2}\,(W(p) + W(-p)) \quad \text{and} \quad
\pi^- W(p) = \tfrac{1}{2}\,(W(p) - W(-p)).
\end{split}
\end{align*}
Then $V$ and $\fourier W$ have the same law.
\end{lemma}

The lemma is proved by the following observations. As both random fields are
centered Gaussian, one only has to show that they have the same covariance. First we
calculate
\begin{align*}
\begin{split}
\expct(\pi^+W(p)& \pi^- W(q))\\
& = \tfrac{1}{4} \,\expct((W(p) + W(-p))(W(q) - W(-q)))\\
& = \tfrac{1}{4} \left( \expct(W(p) W(q)) + \expct(W(-p)W(q)) - \expct(W(p)W(-q)) - \expct(W(-p)W(-q))\right)\\
& = \tfrac{1}{4} \left( \delta(p-q) + \delta(p+q) - \delta(p+q) -\delta(p-q) \right)
 = 0.
\end{split}
\end{align*}
Therefore $\pi^+W$ and $\pi^- W$ are uncorrelated, and hence independent. Then this
observation yields
\begin{align*}
\begin{split}
\expct(V(p)V(q))
& = \expct(\pi^+ W(p) \pi^+ W(q)) - \expct(\pi^- W(p) \pi^- W(q))\\
& = \tfrac{1}{4} (2 \delta(p-q) + 2 \delta(p+q) - (2 \delta(p-q) - 2 \delta(p+q)))
 = \delta(p+q)
\end{split}
\end{align*}
and similarly
\begin{align*}
\expct(V(p) \overline{V(q)}) = \delta(p-q),
\end{align*}
and $V = \fourier W$ in distribution as claimed in Lemma~\ref{lemma_FW}.\qed

By Lemma~\ref{lemma_FW} a random field $\psi$ given by
\begin{equation}\label{eqn_second_GRF_alg}
\psi(x) := ( \ifourier \gamma^{1/2} (\pi^+ W + i \pi^-W))(x)
\end{equation}
is equal in distribution to $\varphi$, and in particular it has covariance~$C$, too.
So Equation~\eqref{eqn_second_GRF_alg} yields a second algorithm for the
construction of random fields with covariance~$C$ which is also presented in
Section~\ref{section_algorithm}. For the following we remark that $V$ as defined
above satisfies  $V(-x) = \overline{V(x)}$ for all $x\in\R^d$. Furthermore, by the
previous calculations, $\Re V(x)$ and $\Im V(x)$ are independent centered Gaussian
random variables.

In the remainder of this section we shall consider the problem of discretizing the
random fields above. Let $A$ be a rectangular region in $\R^d$ defined as the
Cartesian product of $d$ closed intervals $J_i = [a_i,b_i]$, $i=1$, \dots, $d$. For
each $i=1$, \dots, $d$, choose $N_i\in \N$ points $y^i_0 < y^i_1 < \ldots <
y^i_{N_i-1}$ in $[a_i, b_i)$ defining a partition of $J_i$. From now on fix $N =
(N_1,\ldots, N_d)$, and suppose that every $N_i$ is even. We denote by $\cK^N$ the
set of all vectors $K = (k_1,\ldots, k_d)$ with $k_i = 0,\ldots, N_i - 1$. Define
$x_K$ to be the point in $A$ whose $i$--th coordinate $(x_K)_i$ is equal to
$y^i_{k_i}$. Thus the $x_K$, $K\in\cK^N$, form a rectangular grid in $A$. For
$K\in\cK^N$, define the cell $\Delta_K$ in $A$ by
\begin{equation*}
 \Delta_K = \bigtimes_{i=1}^d \, [x_K, x_{K+e_i}),
\end{equation*}
where $e_i = (\delta_{ij}, j=1, \dotsc, d)$. Thus the set $\{\Delta_K,\,K\in\cK^N\}$
of cells forms a partition of $A' = [a_1,b_1)\times\ldots\times[a_d,b_d)$.

As above, we let $W$ denote a white noise random field on $\R^d$, which defines a
discretized white noise random field $W^N$ by
\begin{equation*}
 W^N = (W_K^N, K\in\cK^N)
\end{equation*}
with
\begin{equation*}
 W_K^N = \size{\Delta_K}^{-1}\int_{\Delta_K} W(x)\,dx.
\end{equation*}
Hence the family $\{W_K^N\}$ is an independent family of centered Gaussian random
variables, $W_K^N$ having variance $\size{\Delta_K}^{-1}$. We may also view $W^N$ as
a discrete white noise random field on $A'$ by considering it as being constantly
equal to $W_K^N$ on the cell $\Delta_K$.

Now we impose periodic boundary conditions at the boundary of $A$ as will be
automatically implemented by use of the FFT in our algorithms in the next section.
That means, in the sequel we work on the torus defined by  $A$.  Therefore all
calculations will be done modulo the side lengths of the rectangle $A$, or modulo
the vector $N$, respectively. In particular, the following discrete analogues of
$\Re V$, $\Im V$ in Lemma~\ref{lemma_FW} are well-defined for all $K\in\cK^N$
\begin{equation*}
	R^N_K = \frac{1}{2}\bigl(W^N_K + W^N_{-K}\bigr),\qquad
	I^N_K = \frac{1}{2}\bigl(W^N_K - W^N_{-K}\bigr).
\end{equation*}
Set $V^N = (V^N_K,\,K\in\cK^N)$ with
\begin{equation*}
	V^N_K = R^N_K + i I^N_K.
\end{equation*}
Observe that we have $I^N_{-K}=-I^N_K$ for all $K\in\cK^N$, and therefore
\begin{equation}	\label{conj}
	V^N_{-K} = \overline{V^N_K}.
\end{equation}
In particular, we find that for all $K\in\cK^N$ which are such that
$K=-K$ (mod $N$), $V^N_K$ is real. Figure~\ref{grid} shows an equidistant
grid with $N=(4,4)$ on the unit square of $\R^2$. The solid points are
the grid points, and those which are blue have the property $K=-K$
(mod $(4,4)$). Consequently the random field is real at the blue grid points.
The red points form an example of two points associated via $K$ and $-K$
with each other: $K=(2,1)=-(2,3)$ (mod $(4,4)$). Thus the value of
the random field at one of the red points is equal to the complex conjugate
at the other.
\begin{figure}[ht]
\begin{center}
\includegraphics[scale=.7]{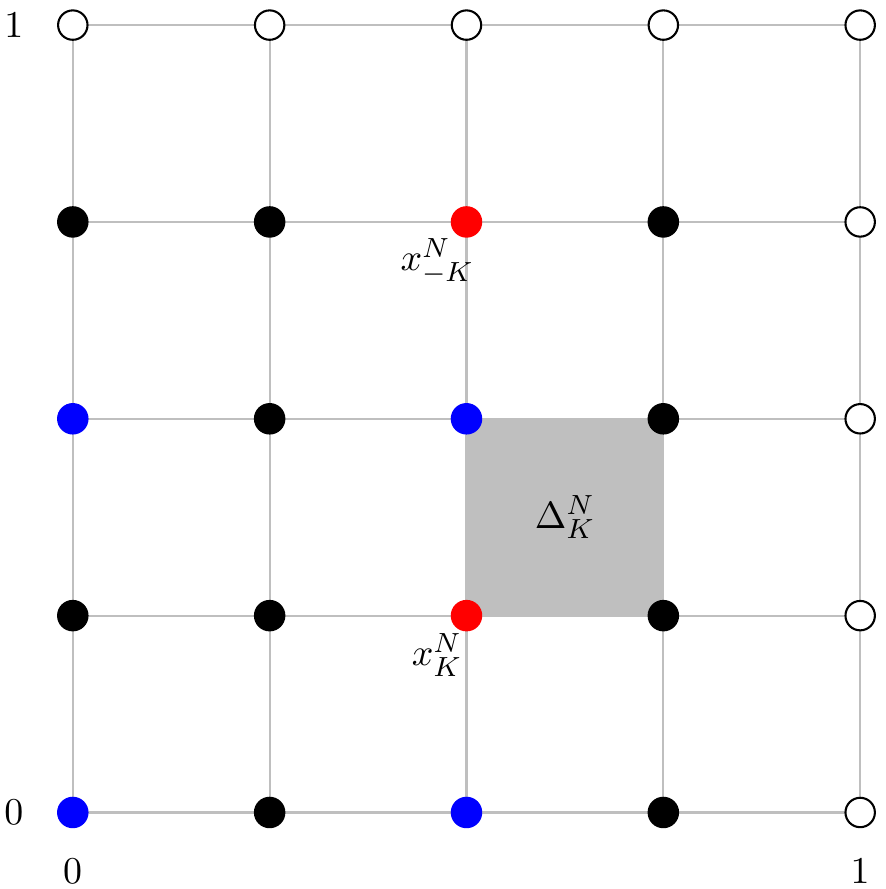}
\caption{A $4\times 4$ grid on the unit square in $\R^2$}\label{grid}
\end{center}
\end{figure}

For $K\in\cK^N$ with $K \ne -K$, the independence of $W^N_K$ and $W^N_{-K}$
entails
\begin{align*}
 \expct\bigl((R^N_K)^2\bigr)
 	& = \expct\Bigr( \bigr(\tfrac{1}{2} \, ( W_K^N + W_{-K}^N)\bigl)^2\Bigl)
 	= \tfrac{1}{4} \, \Bigl( \expct\bigl( (W^N_K)^2\bigr) + 2 \, \expct
		(W_K^N \, W_{-K}^N) + \expct\bigl( (W^N_{-K})^2\bigr) \Bigr)\\
	& = \tfrac{1}{4} \, \bigl( \size{\Delta_K^N}^{-1} +
		\size{\Delta_{-K}^N}^{-1} \bigr),
\end{align*}
and similarly
\begin{equation*}
 \expct\bigl((I^N_K)^2\bigr) =
		\tfrac{1}{4} \, \bigl( \size{\Delta_K^N}^{-1} +
		\size{\Delta_{-K}^N}^{-1} \bigr).
\end{equation*}
On the other hand for $K\in\cK^N$ with $K=-K$ we get
\begin{align*}
\expct\bigl((R^N_K)^2\bigr) = \size{\Delta^N_K}^{-1},
\end{align*}
and --- as already mentioned above --- $I^N_K=0$.

It is straightforward to check that for every $K\in\cK^N$ there is exactly one
vector $L\in\cK^N$ so that $L=-K$ (mod $N$). Let $\cL^N$ denote any maximal subset
of $\cK^N$, so that the random variables $V^N_K$, $K\in\cL^N$ are independent. Thus
we only have to generate the random variables $R^N_K$ and $I^N_K$ for $K\in\cL^N$,
and the values of the random field $V$ at the other indices are determined from
these by complex conjugation. In the appendix it is shown how one can construct one
such set $\cL^N$ for an arbitrary dimension $d\in\N$. Here we only remark that for
$d=2$, a possible choice for $\cL^N$, $N=(N_1,N_2)$, is given by
\begin{equation*}
\begin{split}
	\cL^{(N_1,N_2)}
		= \{k_1=0,&\dotsc,M_1,\,k_2=0,\dotsc,M_2\}\\
			&\uplus \{k_1=1,\dotsc,M_1-1,\,k_2=M_2+1,\dotsc,2M_2-1\},
\end{split}
\end{equation*}
where we have set $M_i=N_i/2$, $i=1$,~$2$. In Figure~\ref{grid_2} this set is shown
for an equipartition of the unit square with an $8\times 8$ grid as those
corresponding to the magenta points. The values of the random field at the green and
black points is determined from those at the magenta points via
Relation~\eqref{conj}.
\begin{figure}[ht]
\begin{center}
\includegraphics[scale=.7]{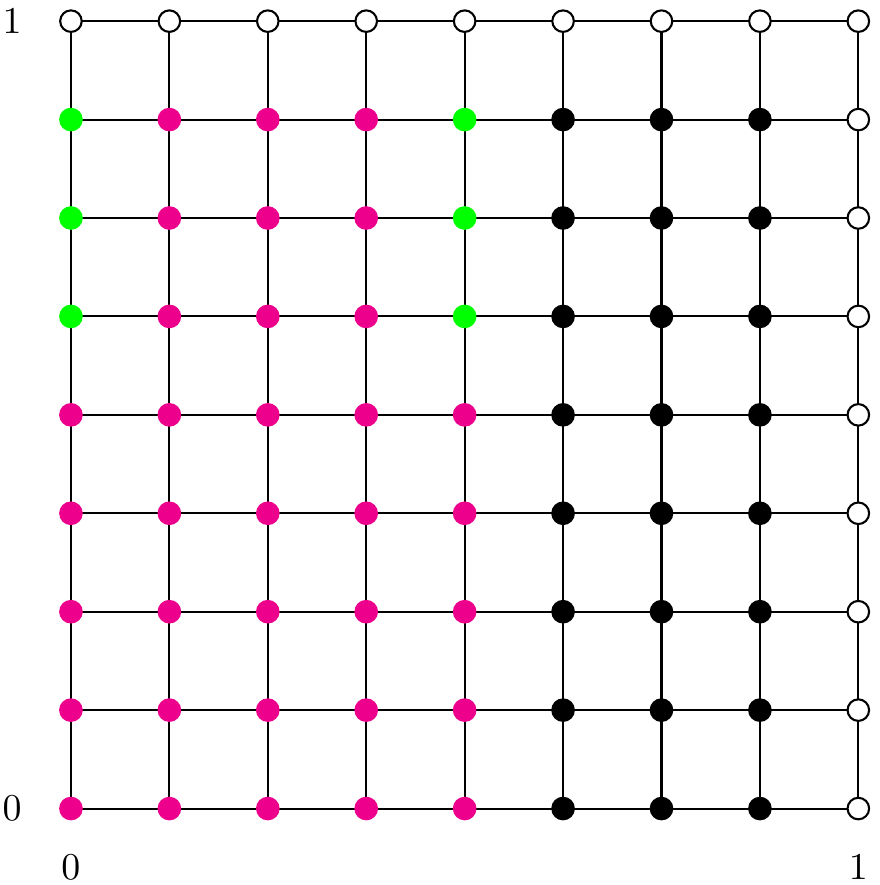}
\caption{$\cL^N$ for an $8\times 8$ grid given by the magenta points}\label{grid_2}
\end{center}
\end{figure}

For simplicity we assume from now on equidistant partitions in each direction, and
let $\size{\Delta^N}$ denote the volume of the cells. Then the calculations above
can be summarized by the following rules:
\begin{enumerate} \renewcommand{\labelenumi}{\alph{enumi})}
	\item In the cell $\Delta^N_K$, $K\in\cK$, the value of the discrete
			random field to be simulated is given by
			\begin{equation}	\label{drf}
				{\rm DFT}^{-1} \gamma(p_K)^{1/2}
				  \bigl(R^N_K + i I^N_K\bigr)|A|^{-1},
			\end{equation}
		where ${\rm DFT}$ denotes the discrete Fourier transform.
	\item the random variables $R^N_K$, $I^N_K$, $K\in\cL^N$, are independent
	      	of each other;		
	\item if $K\in\cL^N$ is such that $K=-K$ (mod $N$), then $R^N_K$ has the law
			$\cN(0,|\Delta^N|^{-1})$, $I^N_K=0$;
	\item if $K\in\cL^N$ is such that $K\ne-K$ (mod $N$), then $R^N_K$
			and $I^N_K$ are distributed with the law
			$\cN(0,(2|\Delta^N|)^{-1})$,
			and $I^N_{-K} = -I^N_K$;
	\item for $K \in \cK$, the points~$p_K$ are chosen such that
			$(p_K)_i = i \cdot (b_i - a_i)$, for $i = 1,2$.
\end{enumerate}
The corresponding algorithm implementing these rules in an efficient way can be
found in the next section.

We close this section with the following

\begin{remark}\label{rem_add_rand_numb}
For higher dimensions $d$ than $d=2$, the implementation of $\cL^N$ becomes a rather
tedious task, see the construction of~$\cL^N$ in the appendix. Consider again the
case $d=2$ in figure~\ref{grid_2}. If we add to $\cL^N$ the set of green points,
i.e., some of the points at the boundary of the upper left subcube, then this set
can be enumerated in a very simple way as $k_1=0,\dotsc,M_1$, $k_2=0,\dotsc,2M_2-1$.
After generation of the random variables, the values at the green points have to be
overwritten according to Relations~\eqref{conj}. Since the random variables which
are superfluously generated are those at the \emph{boundary} of the upper left cube,
their number is negligible as compared to those in the interior of that cube. A
similar consideration holds in the case of general $d\in\N$: If one generates
additionally to the random variables indexed by $\cL^N$ those at the boundaries of
the subcubes, one obtains a simple code at negligibly higher computing costs. After
generation of all values, those which have been superfluously generated have to be
overwritten with the correct values according to Equation~\eqref{conj}.
\end{remark}

\section{Algorithms}\label{section_algorithm}

The following algorithms are based on the calculations of the previous section. They
allow a fast generation of $d$-dimensional stationary Gaussian random fields with
given covariance using the advantage of fast Fourier transforms instead of directly
calculating computational expensive convolutions. The boundary conditions
implemented here are periodic. For Neumann or Dirichlet boundary conditions, fast
cosine and sine transformations have to be used.

The first algorithm is an implementation of Equation~(\ref{eqn_first_GRF_alg}).

\begin{algorithm}\label{algorithm_two_FFT}
$ $ \newline
\emph{Remarks:} \nopagebreak
\begin{enumerate}
\item The functions $\fft$ and $\ifft$ include all necessary rescaling depending
on the used $\fft$ algorithm and the integers $N_i$.
\item $A$ is a $d$-dimensional complex-valued array, $B$ is real-valued.
\item $x_{k_1\cdots k_d}$ denotes the grid point corresponding to the integers
$(k_1,\ldots,k_d)$. The grid points are distributed equidistantly in each
direction, i.e.\ the distance of two arbitrary neighbor grid points in direction
$e_i$ is given by a constant $\Delta x_i$.
\item The points $p_{k_1\cdots k_d}$ in the Fourier domain are given by
$(p_{k_1\cdots k_d})_i = (k_i-N_i/2)/l_i$ for $i=1,\ldots, d$.
\end{enumerate}
\emph{Input:} \begin{enumerate}
		\item $d$-dimensional rectangular region $D$, where
			$l_1,\ldots,l_d$ is the length of the edges,
		\item $N_1,\ldots,N_d$ number of discretization points in each
direction, all even,
		\item $\gamma^{1/2}$ a symmetric, positive function on $\R^d$,
		\item $R()$ a function that generates independent
		      $\cN(0,|\Delta^N|^{-1})$-distributed random numbers.
              \end{enumerate}
\emph{Output:} GRF $B$ on $D$, where the covariance is given by the Fourier
transform of $\gamma$.

\vspace*{\baselineskip}
\begin{algorithmic}
\FOR{$k_i = 0, \ldots, N_i - 1, \, i=1,\ldots,d$}
\STATE $B(k_1,\ldots,k_d) \leftarrow R();$
\ENDFOR
\STATE $A \leftarrow \fft \, B;$
\FOR{$k_i = 0, \ldots, N_i - 1, \, i=1,\ldots,d$}
\STATE $A(k_1,\ldots,k_d) \leftarrow A(k_1,\ldots,k_d) \cdot
	\gamma(p_{k_1\cdots k_d})^{1/2} \cdot |D|^{-1};$
\ENDFOR
\STATE $B \leftarrow \ifft A;$
\end{algorithmic}
\end{algorithm}

The second algorithm is based on Equation~(\ref{drf}) with
Remark~\ref{rem_add_rand_numb}. For the direct construction of the Fourier transform
of the white noise field, a random field with even real and odd imaginary part has
to be constructed. At all grid points indexed by $K\in\cL^N$ with $K=-K$
$\text{mod}(N)$ the random field is real-valued as the random field is equal to its
complex conjugate at these points. At the grid points where the random field is
real-valued it is scaled with a factor one instead of $2^{-1/2}$. If we use a
$d$-dimensional DFT on a grid with $N_i$ discrete points in direction $e_i$, the
following algorithm generates a GRF with approximate covariance $C$. In order to
keep the error small, one has to take care that the correlation length defined by
the covariance is small in comparison to the size of the rectangular region.
\goodbreak

\begin{algorithm} \label{algorithm_one_FFT}
$ $\newline
\emph{Remarks:}
\begin{enumerate}
\item All calculations have to be done modulo $N_i$ in the $i$--th direction.
\item The function $\ifft$ includes all necessary rescaling depending on the used $\fft$
      algorithm and the integers $N_i$.
\item $A$ is a $d$-dimensional complex-valued array, $B$ is real-valued.
\item $x_{k_1\cdots k_d}$ denotes the grid point corresponding to the integers
		$(k_1,\ldots,k_d)$. The grid points are distributed equidistantly
         in each direction, i.e.\ the distance of two arbitrary neighboring grid
         points in direction $e_i$ is given by a constant $\Delta x_i$.
\item The points $p_{k_1\cdots k_d}$ in the Fourier domain are given by
$(p_{k_1\cdots k_d})_i = (k_i-N_i/2)/l_i$ for $i=1,\ldots, d$.
\end{enumerate}
\emph{Input:}
 \begin{enumerate}
		\item $d$-dimensional rectangular region $D$ with
			$l_1,\ldots,l_d$ length of the edges,
		\item $N_1,\ldots,N_d$ number of discretization points in each direction, all even,
		\item $\gamma^{1/2}$ a symmetric, positive function on $\R^d$,
		\item $R()$ a function that generates independent
		      $\cN(0,|\Delta^N|^{-1})$-distributed random numbers.
              \end{enumerate}
\emph{Output:} GRF $B$ on $D$, where the covariance is given by the Fourier transform of $\gamma$.

\vspace*{\baselineskip}
\begin{algorithmic}
\FOR{$k_i = 0, \ldots, N_i-1, \, i=1,\ldots,d-1$, $k_d = 0, \ldots, N_d/2$}
\IF{$k_i \in \set{0,N_i/2}$, for all $i=1,\ldots,d$}
\STATE $\Re A(k_1,\ldots,k_d) \leftarrow R()\cdot
	  \gamma(p_{k_1\cdots k_d})^{1/2} \cdot |D|^{-1};$
\STATE $\Im A(k_1,\ldots,k_d) \leftarrow 0;$
\ELSE
\STATE $\Re A(k_1,\ldots,k_d) \leftarrow 2^{-1/2} R()\cdot
	    \gamma(p_{k_1\cdots k_d})^{1/2} \cdot |D|^{-1};$
\STATE $\Im A(k_1,\ldots,k_d) \leftarrow 2^{-1/2} R()\cdot
	    \gamma(p_{k_1\cdots k_d})^{1/2} \cdot |D|^{-1};$
\STATE $\Re A(N_1 - k_1,\ldots,N_d - k_d) \leftarrow \Re A(k_1,\ldots,k_d);$
\STATE $\Im A(N_1 - k_1,\ldots,N_d - k_d) \leftarrow - \Im A(k_1,\ldots,k_d);$
\ENDIF
\ENDFOR
\STATE $B \leftarrow \ifft A;$
\end{algorithmic}
\end{algorithm}

Next, we give an algorithm that generates just the necessary amount of random
numbers. For that algorithm, we observe that the two hyperplanes induced by $k_i =
0, \ldots, N_i-1$ for $i = 1, \ldots, d-1$ and $k_d = 0$ or $k_d = N_d/2$ have the
same structure as the corresponding torus in $d-1$~dimensions. So we use a recursive
implementation that solves the problem of setting the values of the random field at
the grid points by setting all values except those in the two hyperplanes, and then
restarts the algorithm for these hyperplanes in $d-1$~dimensions.

\begin{algorithm}
\label{algorithm_one_FFT_recursive}
\ \newline
\emph{Remarks:}
\begin{enumerate}
\item All calculations have to be done modulo $N_i$ in the $i$-th direction.
\item The function $\ifft$ includes all necessary rescaling depending on the
used $\fft$ algorithm and the integers $N_i$.
\item $A$ is a $d$-dimensional complex-valued array, $B$ is real-valued.
\item $x_{k_1\cdots k_d}$ denotes the grid point corresponding to the integers
$(k_1,\ldots,k_d)$. The grid points are distributed equidistantly in each
direction, i.e.\ the distance of two arbitrary neighbor grid points in direction
$e_i$ is given by a constant $\Delta x_i$.
\item The points $p_{k_1\cdots k_d}$ in the Fourier domain are given by
$(p_{k_1\cdots k_d})_i = (k_i-N_i/2)/l_i$ for $i=1,\ldots, d$.
\end{enumerate}
\emph{Input:}
 \begin{enumerate}
		\item $d$-dimensional rectangular region $D$ with
			$l_1,\ldots,l_d$ length of the edges,
		\item $N_1,\ldots,N_d$ number of discretization points in each
direction, all even,
		\item $\gamma^{1/2}$ a symmetric, positive function on $\R^d$,
		\item $R()$ a function that generates independent
$\cN(0,|\Delta^N|^{-1})$-distributed random numbers.
              \end{enumerate}
\emph{Output:} GRF $B$ on $D$, where the covariance is given by the Fourier
transform of $\gamma$.

\begin{algorithmic}
 \STATE \textbf{function} \emph{create\_random\_field}($D,N_1,\ldots,N_d$)

  \STATE $A \leftarrow $ complex\_array$(N_1, \ldots, N_d)$;
  \STATE \emph{set\_array}($\{N_1,\ldots,N_d\}$, $\{ \; \}$, $A$);

\STATE $B \leftarrow \ifft A;$

  \STATE \textbf{end function}
\end{algorithmic}

\vspace*{0.5\baselineskip}

\begin{algorithmic}
\STATE \textbf{function} \emph{set\_array}($\{N_1,\ldots,N_j\}$, $\{ k_{j+1},
\ldots, k_d\}$, $A$)
\FOR{$k_i = 0, \ldots, N_i, \, i=1,\ldots,j-1$, $k_j = 1, \ldots,
N_j/2-1$}
\STATE $\Re A(k_1,\ldots,k_d) \leftarrow  2^{-1/2} R()\cdot
\gamma(p_{k_1\cdots k_d})^{1/2} \cdot |D|^{-1};$
\STATE $\Im A(k_1,\ldots,k_d) \leftarrow  2^{-1/2} R()\cdot
\gamma(p_{k_1\cdots k_d})^{1/2} \cdot |D|^{-1};$
\STATE $\Re A(N_1 - k_1,\ldots,N_d - k_d) \leftarrow \Re A(k_1,\ldots,k_d);$
\STATE $\Im A(N_1 - k_1,\ldots,N_d - k_d) \leftarrow - \Im A(k_1,\ldots,k_d);$
\IF{($|\{N_1,\ldots,N_j\}| = 1$)}
\STATE $\Re A(0,\ldots,k_d) \leftarrow R()\cdot \gamma(p_{0 \, k_2\cdots
k_d})^{1/2} \cdot |D|^{-1};$
\STATE $\Im A(0,\ldots,k_d) \leftarrow 0;$
\STATE $\Re A(N_1/2,\ldots,k_d) \leftarrow R()\cdot
\gamma(p_{N_1/2 \, k_2\cdots k_d})^{1/2} \cdot |D|^{-1};$
\STATE $\Im A(N_1/2,\ldots,k_d) \leftarrow 0;$
\ELSE
\STATE \emph{set\_array}($\{N_1,\ldots, N_{j-1}\}$, $\{ 0, k_{j+1}, \ldots, k_d\}$,
$A$);
\STATE \emph{set\_array}($\{N_1,\ldots, N_{j-1}\}$, $\{ N_{j}/2, k_{j+1},
\ldots, k_d\}$, $A$);
\ENDIF
\ENDFOR
\STATE \textbf{end function}
\end{algorithmic}
\end{algorithm}

Finally we analyze the computational costs of the algorithms in this section. Let $N
= N_1 + \cdots + N_d$. Algorithm~\ref{algorithm_two_FFT} sets the random array in
$\Op(N)$, each FFT is done in~$\Op(N \log N)$ and the costs for the multiplication
with $\gamma^{1/2}$ is each for the real and the imaginary part $\Op(N)$. So overall
we have $\Op(N \log N)$. The second algorithm~\ref{algorithm_one_FFT} performs the
generation of the complex random array in~$\Op(N)$ for real and imaginary part and
then one FFT in~$\Op(N \log N)$. The additional generated random numbers cost
$\Op(N_1 + \cdots + N_{d-1})$. Therefore the algorithm also has complexity $\Op(N
\log N)$. Nevertheless, the constants are smaller. In simulations,
Algorithm~\ref{algorithm_two_FFT} took twice as long as
Algorithm~\ref{algorithm_one_FFT}. Finally,
Algorithm~\ref{algorithm_one_FFT_recursive} constructs the complex random field
in~$\Op(N)$ and performs the FFT in~$\Op(N \log N)$, which overall leads to a
complexity of~$\Op(N \log N)$. Therefore all algorithms are in the same complexity
class but in simulations they perform with different speeds due to different
constants. So in terms of complexity, there is no gain if one omits one FFT, but in
practice a factor two in speed could be of advantage.

The algorithms have been implemented and tested for $d = 1,2$. The results are
presented in the subsequent section.

\section{Simulations and Statistical Tests}\label{section_tests}

This section shows the results of an implementation in C++ of the algorithms
presented in Section~\ref{section_algorithm} using FFTW~\cite{Frigo:2005}. The
simulations were done on a rectangle in $\R^2$ and on the interval $[-\pi, \pi]$ in
$\R$. We first consider the two-dimensional implementation. Afterwards we do error
analysis on the interval. It turned out that both algorithms lead to the same
results but that Algorithm~\ref{algorithm_one_FFT} is twice as fast as
Algorithm~\ref{algorithm_two_FFT}.

For the two-dimensional implementation, the covariance was chosen to be
\begin{equation*}
C(x,y) = \integral{\R^d}{e^{-2\pi i \sip{p}{x-y}} \gamma(p)}{p}, \qquad x, y \in
\R^2.
\end{equation*}
with
\begin{equation}
\gamma(p) = (m^{kl} + (p_1^{2k} + p_2^{2k})^l)^{-n},
\end{equation}
where $k,l,n \in \N$, $m \in \R_{>0}$, and $p = (p_1,p_2) \in \R^2$. For $k = l = n
= 1$ this covariance function is used in Euclidean quantum field theory,
cf.~\cite{Glimm:1981}. It is shown in~\cite{Lang:2007} that $C$ has exponential
decay for $\size{x-y} \gg 1$. (This is well-known for $k=l=n=1$,
e.g.~\cite{Glimm:1981}.) Statistical tests were made on the samples by choosing one
grid point $x_0$ and calculating the statistical covariance of the corresponding
random field with respect to the random field at all other grid points. Some results
are shown in Figure~\ref{figure_3d_covariance}.
\begin{figure}
\centering
\subfigure[$m=10,n=1,k=1,l=2$]{
\includegraphics[scale=0.22]{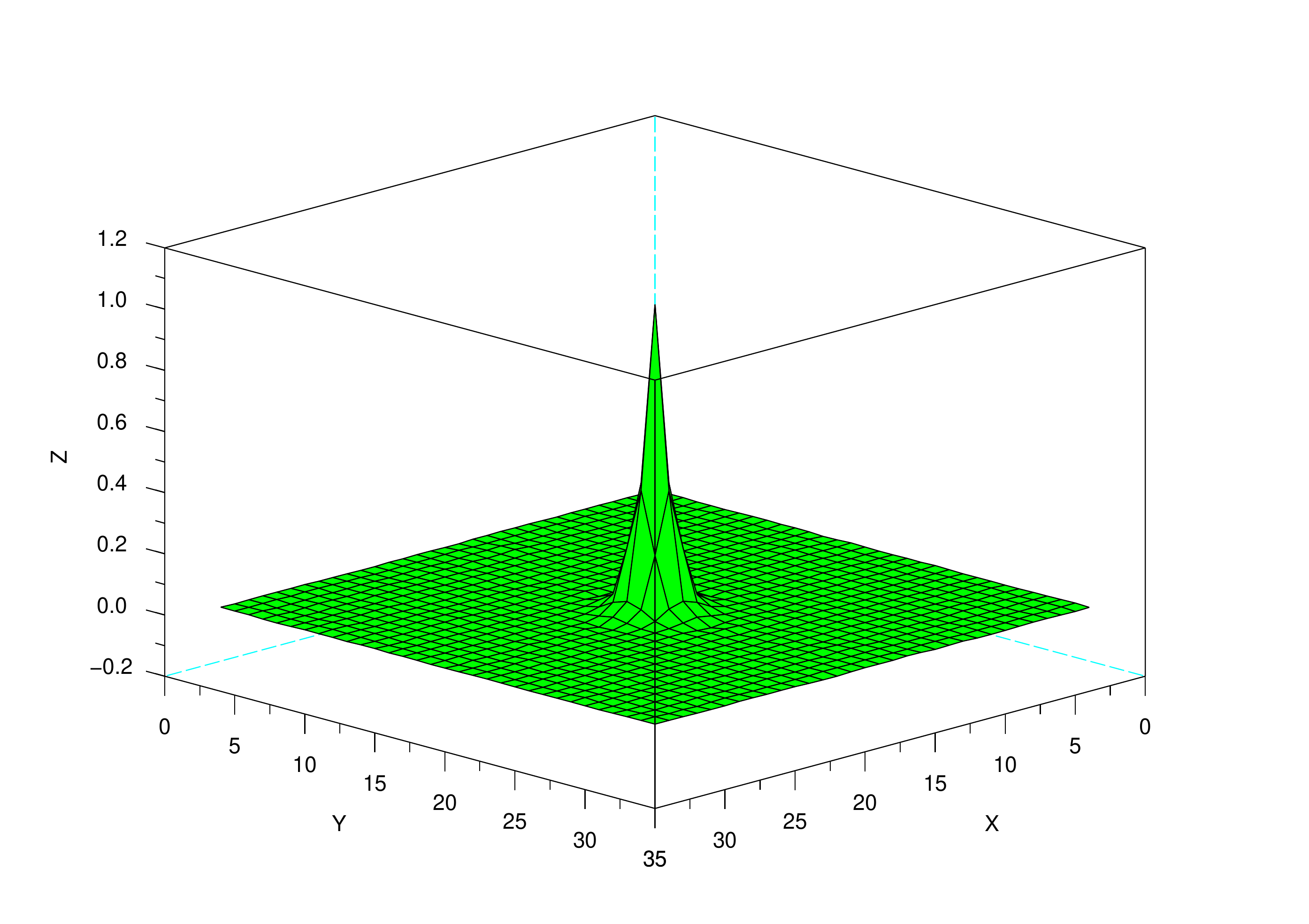}
\label{figure_3d_covariance_green} }
\subfigure[$m=0.75,n=1,k=1,l=2$]{
\includegraphics[scale=0.22]{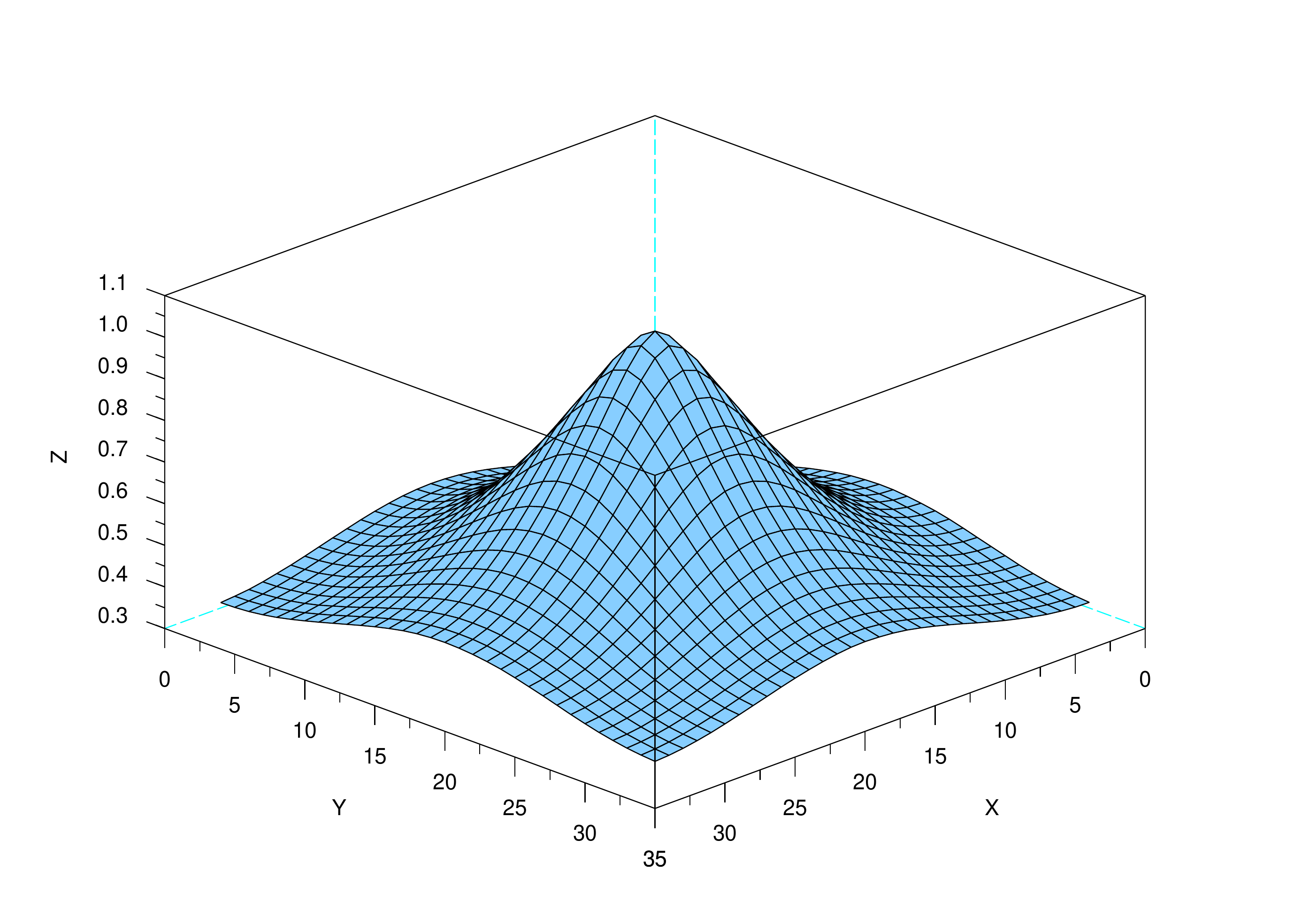}
\label{figure_3d_covariance_blue} }
\subfigure[$m=10,n=3,k=1,l=1$]{
\includegraphics[scale=0.22]{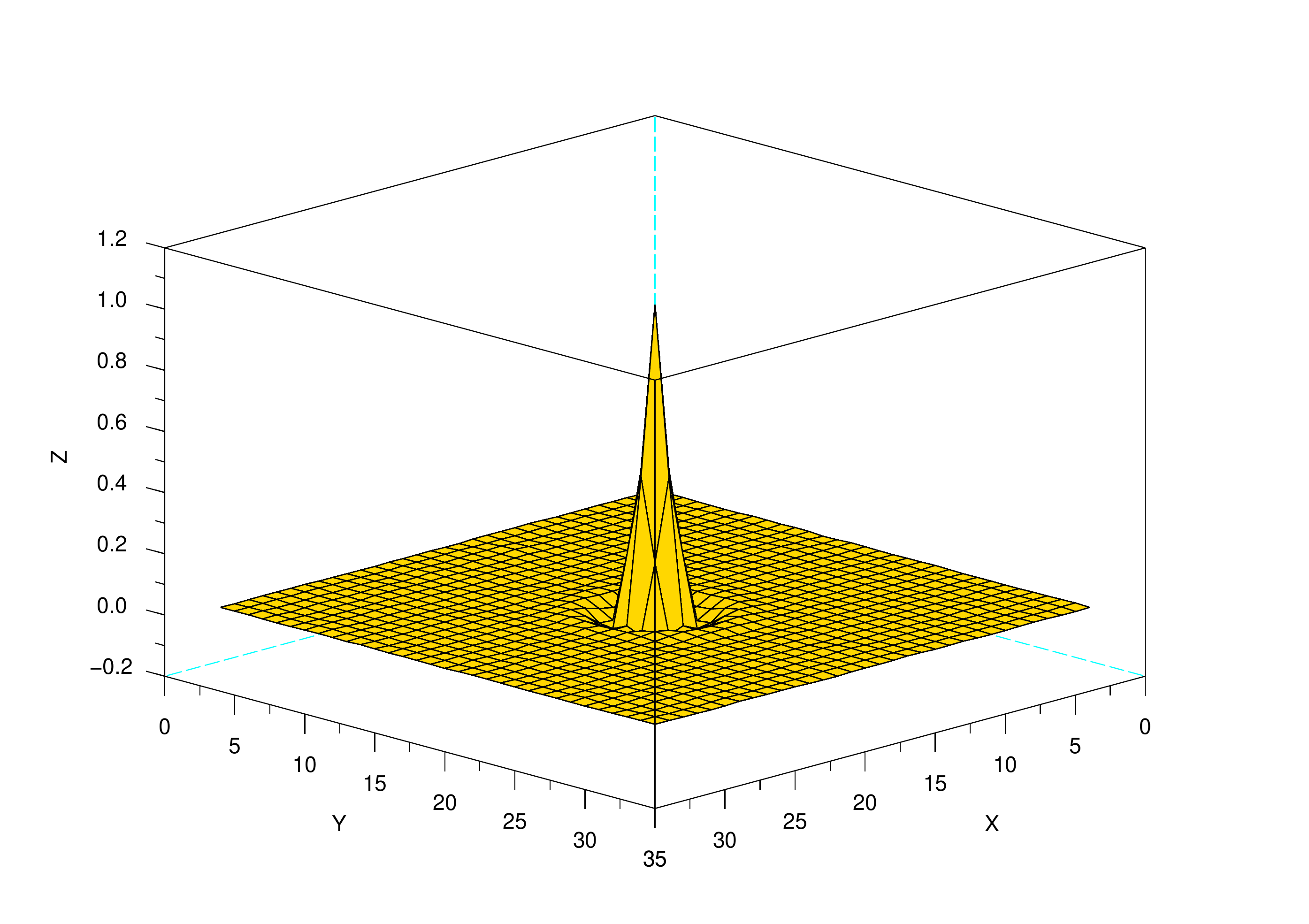}
\label{figure_3d_covariance_yellow} }
\subfigure[$m=5,n=1,k=3,l=1$]{
\includegraphics[scale=0.22]{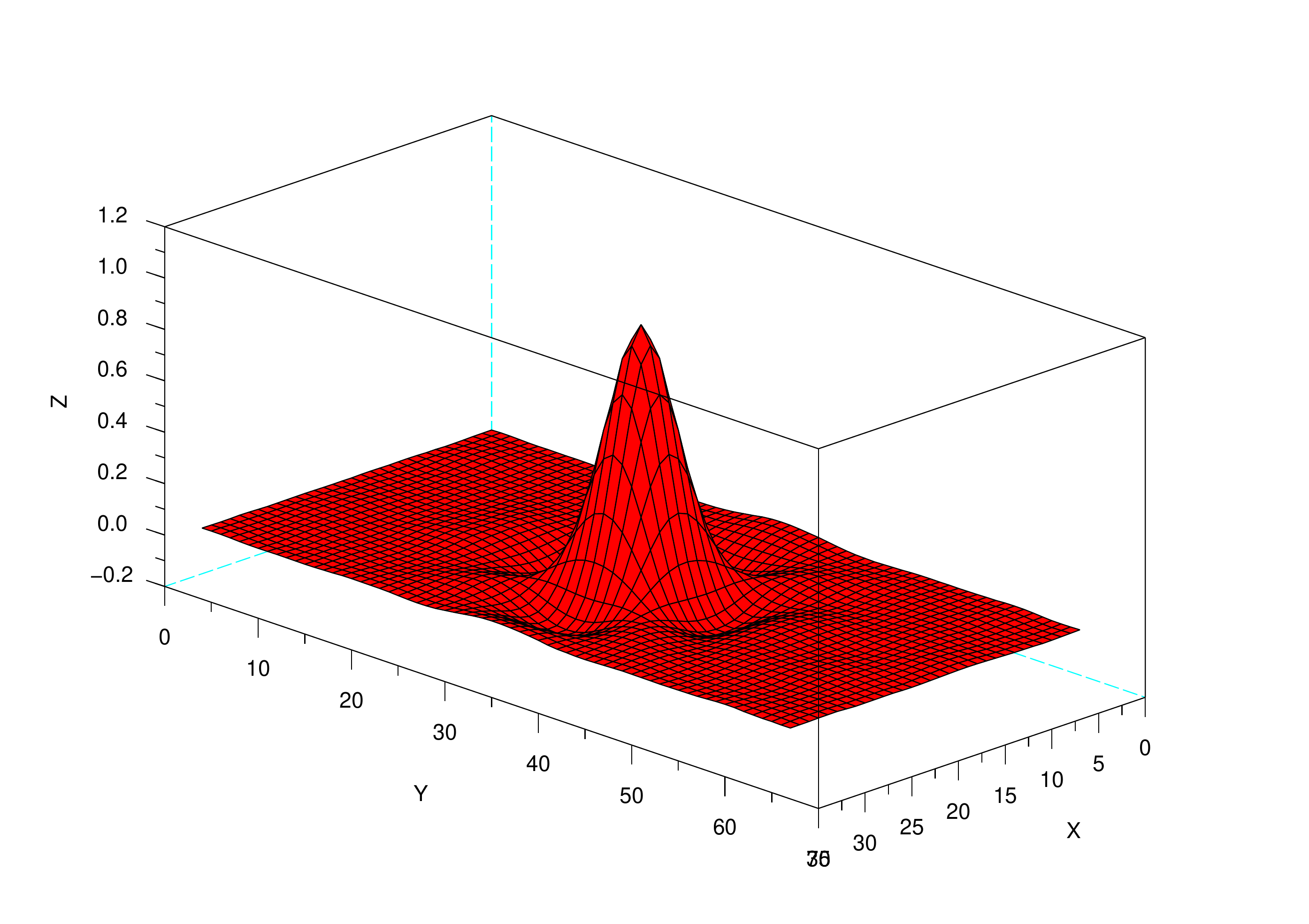}
\label{figure_3d_covariance_red} }
\caption{Simulation results with different covariance
functions.\label{figure_3d_covariance}}
\end{figure}
We remark that the graph in Figure~\ref{figure_3d_covariance_yellow} shows negative
covariance values, however this is the correct behaviour and not an artifact of the
simulation. A physical interpretation of that phenomenon would be some kind of
countermove against the displacement of a point like an elastic band.
Figure~\ref{figure_3d_covariance_red} shows a covariance function that is not
symmetric under rotation. Functions like that might be interesting for applications
with direction dependent correlations.

The simulated random fields are visualized in Figure~\ref{figure_3d_2d_GRF}, \ref{figure_GRF_same_degree},
and \ref{figure_GRF_diff_degree}. Two different possibilities how to visualize the data are presented in
Figure~\ref{figure_3d_2d_GRF}.
\begin{figure}
\centering
\subfigure[3-dimensional view]{
\includegraphics[scale=0.25]{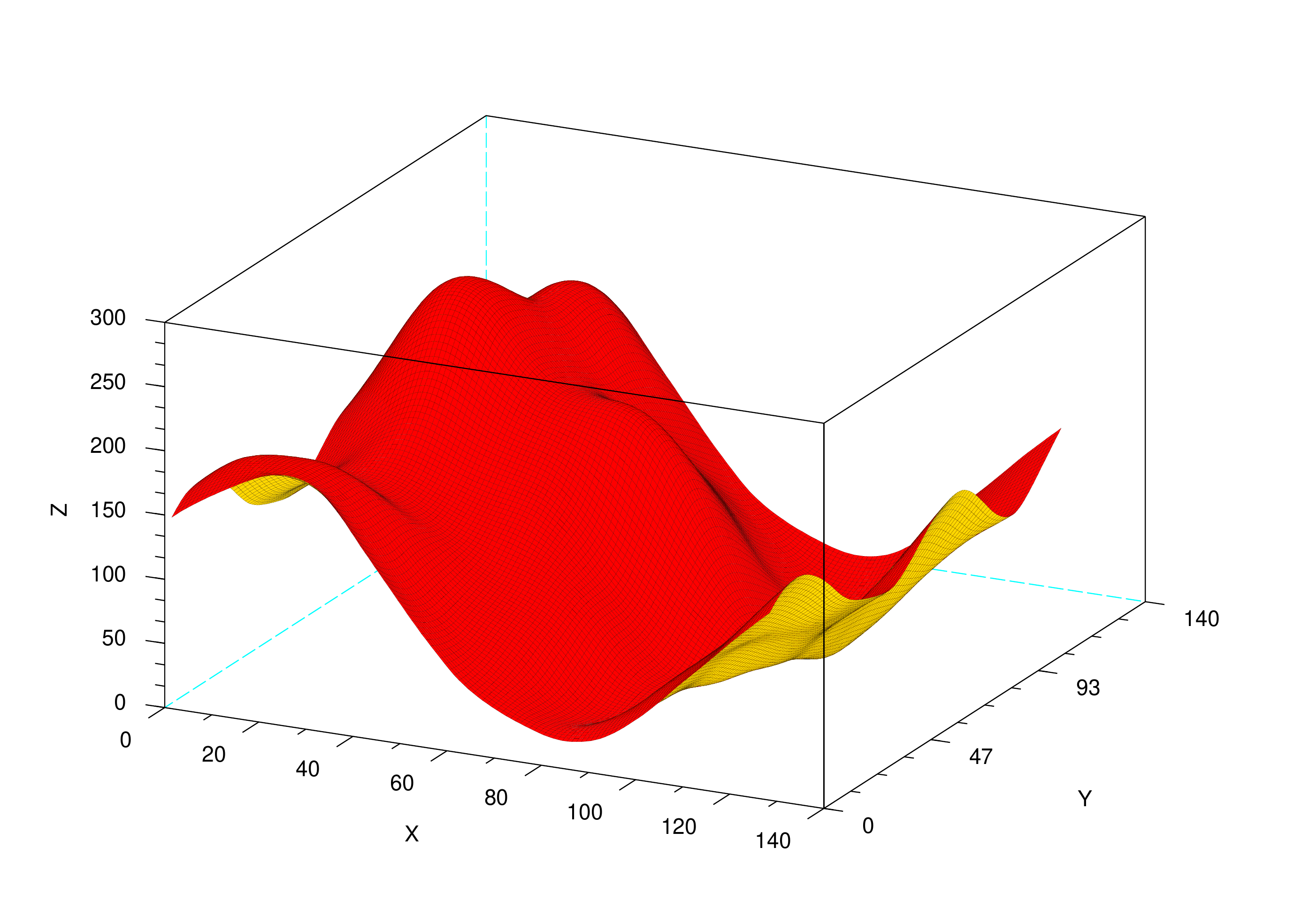}}
\subfigure[2-dimensional implicit function]{
\includegraphics[scale=0.32,origin=b]{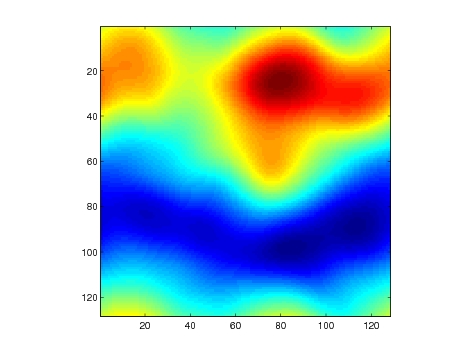}}
\caption{Two different ways of visualizing the same GRF.\label{figure_3d_2d_GRF}}
\end{figure}
On the left, the random numbers are plotted as the graph of a function from $\R^2$
to $\R$. The figure on the right-hand side presents the random numbers in the form
of colors where blue stands for small numbers and red for the large ones.
Figure~\ref{figure_GRF_same_degree}
\begin{figure}
\centering
\subfigure[$\gamma(p) = (1 + \size{p}^2)^{-1}$]{
\includegraphics[scale=0.3]{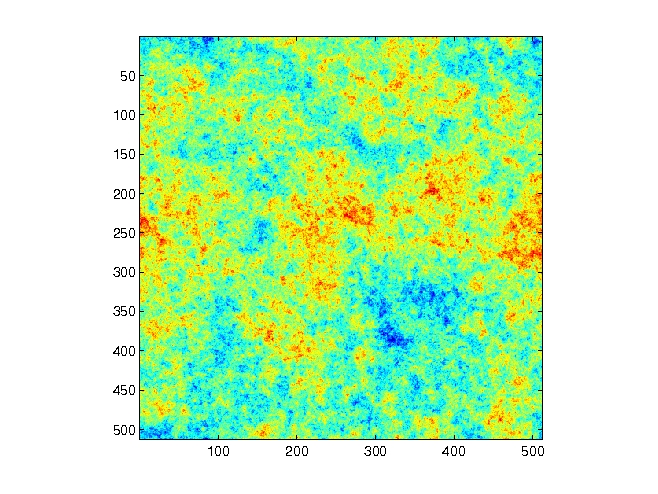}\label{figure_GRF_Glimm}}
\subfigure[$\gamma(p) = (1 + \size{p}^2)^{-2}$]{
\includegraphics[scale=0.3]{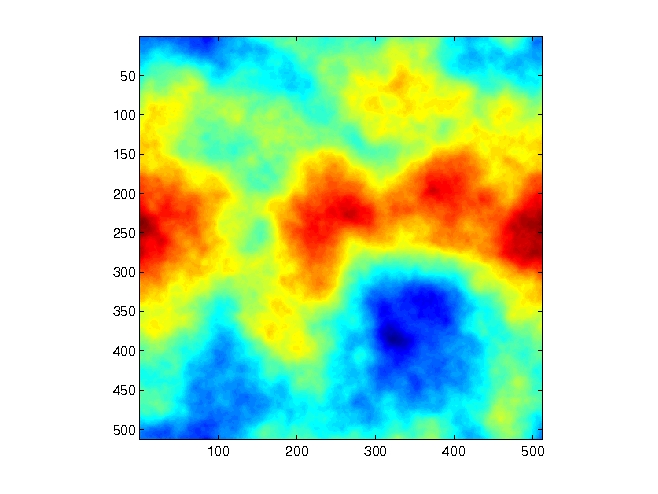}}
\subfigure[$\gamma(p) = (1 + \size{p}^4)^{-1}$]{
\includegraphics[scale=0.3]{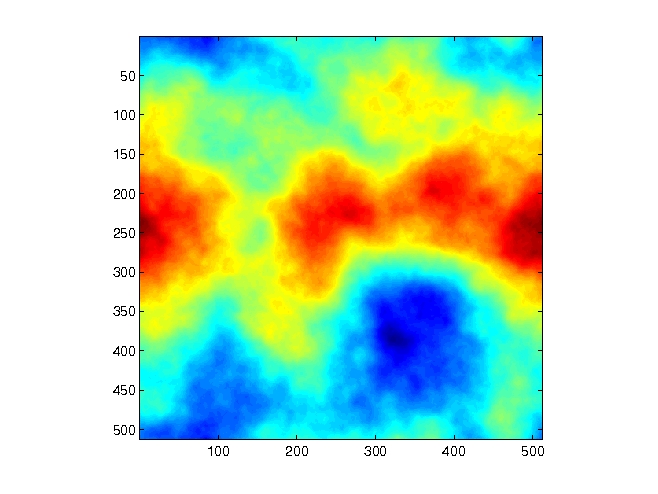}}
\subfigure[$\gamma(p) = (1 + p_1^4 + p_2^4)^{-1}$]{
\includegraphics[scale=0.3]{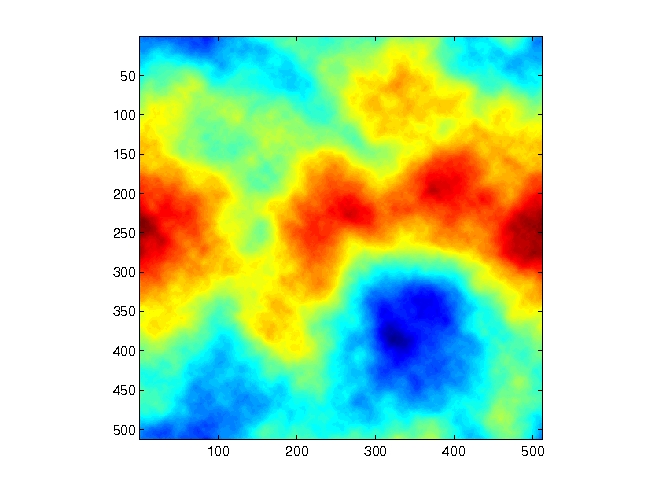}}
\caption{Images of size $512\times512$ with the same white noise.\label{figure_GRF_same_degree}}
\end{figure}
uses this hot color map to compare correlated random fields that use the same white
noise samples for $W$ but different covariance functions. There are differences
between all of the pictures, but it attracts attention that there are hardly any
differences between the three pictures where $\gamma$ is of degree $-4$, while
Picture~\ref{figure_GRF_Glimm} where $\gamma$ is of degree~$-2$ is completely
different which is due to infinite variance of the underlying random field. This
phenomenon can be observed for all tested other degrees as well.
Figure~\ref{figure_GRF_diff_degree}
\begin{figure}
\centering
\subfigure[$\gamma(p) = (1 + \size{p}^2)^{-1}$]{
\includegraphics[scale=0.29]{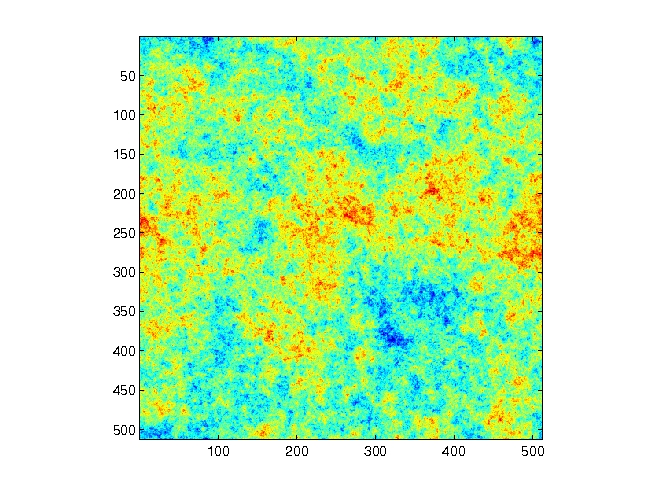}}
\subfigure[$\gamma(p) = (1 + \size{p}^2)^{-2}$]{
\includegraphics[scale=0.29]{case2_2.jpg}}
\subfigure[$\gamma(p) = (1 + \size{p}^2)^{-3}$]{
\includegraphics[scale=0.29]{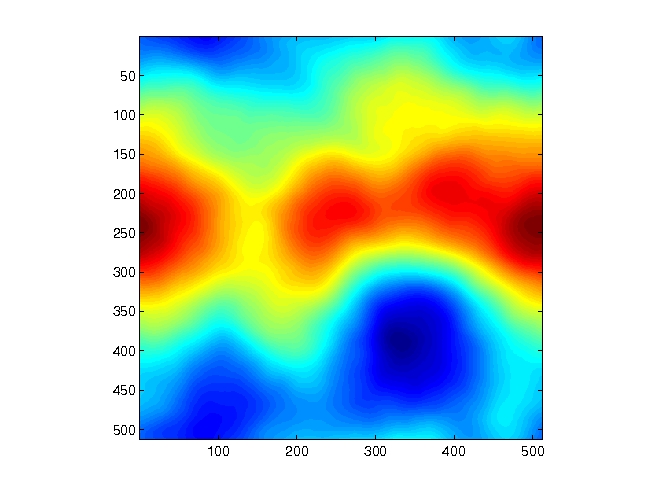}}
\subfigure[$\gamma(p) = (1 + \size{p}^2)^{-4}$]{
\includegraphics[scale=0.29]{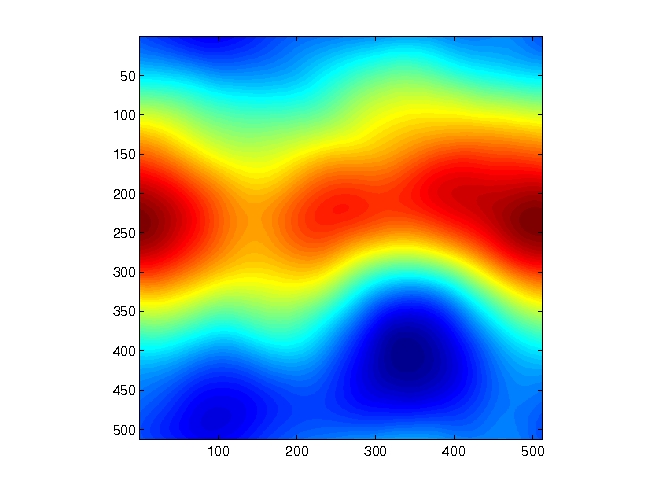}}
\subfigure[$\gamma(p) = (1 + \size{p}^2)^{-5}$]{
\includegraphics[scale=0.29]{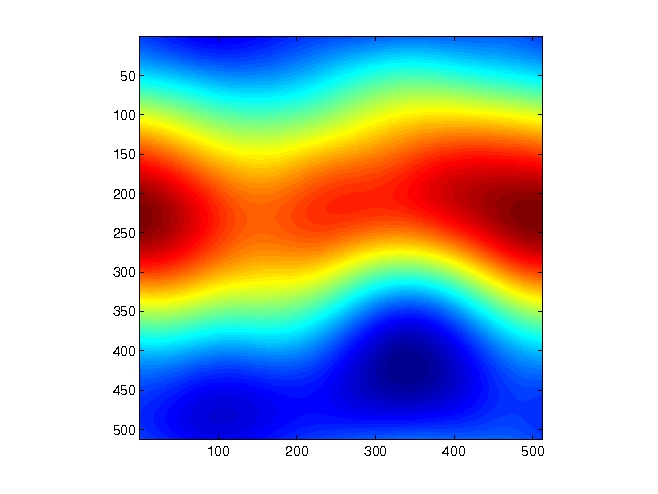}}
\subfigure[$\gamma(p) = (1 + \size{p}^2)^{-6}$]{
\includegraphics[scale=0.29]{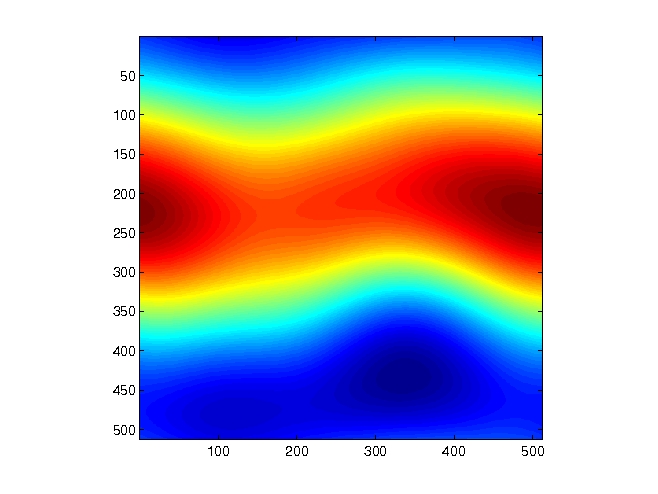}}
\caption{Images of size $512\times512$ with the same white noise.\label{figure_GRF_diff_degree}}
\end{figure}
shows that using the same white noise and the same type of covariance function with
different degrees of $\gamma$ leads to similar samples of the random fields but the
higher the degree of the polynomial, the smoother the resulting noise as expected.
Setting $m$ to something larger than~$1$ will result in faster decreasing
correlations which can be seen in Figure~\ref{figure_3d_covariance}.

We analyze the convergence of the statistical covariance on the interval
$[-\pi,\pi] \subset \R$ with theoretical covariance~$C$ defined by
\begin{align*}
 \gamma(p) = \frac{32 \cdot 10^6}{\pi} \cdot \frac{1}{(100 + p^2)^4},
\end{align*}
for $p \in \R$. The constant is used to scale the variance to~$1$. Then $C$ can be
computed explicitly, and it is given by
\begin{align*}
 C(x,y) = ( 200/3 |x-y|^3 + 40 |x-y|^2 + 10 |x-y| + 1) \exp(-10 |x-y|),
\end{align*}
for $x,y \in \R$. For the error analysis we simulated $10^8$~samples
$(\eta^n_j, j=1,\ldots, 10^8)$ on the equidistant grids $D^n
= \{x_i^n = i \cdot \pi/2^{n-1}, i = -2^{n-1}, \ldots, 2^{n-1}\}$ for $n = 2,
\ldots, 12$. The covariance was then estimated by
\begin{align*}
\overline{\text{Cov}}(x_i^n) =
\sum_{j = 1}^{10^8}
      \eta^n_j(0) \cdot \eta^n_j(x_i^n)
\end{align*}
for all $i = -2^{n-1}, \ldots, 2^{n-1}$ and $n = 2,\ldots,12$.
Figure~\ref{figure_stat_cov_1d} shows the exact covariance~$C$ and the
statistical results on a grid with $9$, $33$, and $4097$ points.
\begin{figure}
 \centering
\subfigure[theoretical covariance]{
\includegraphics[width=0.48\textwidth]{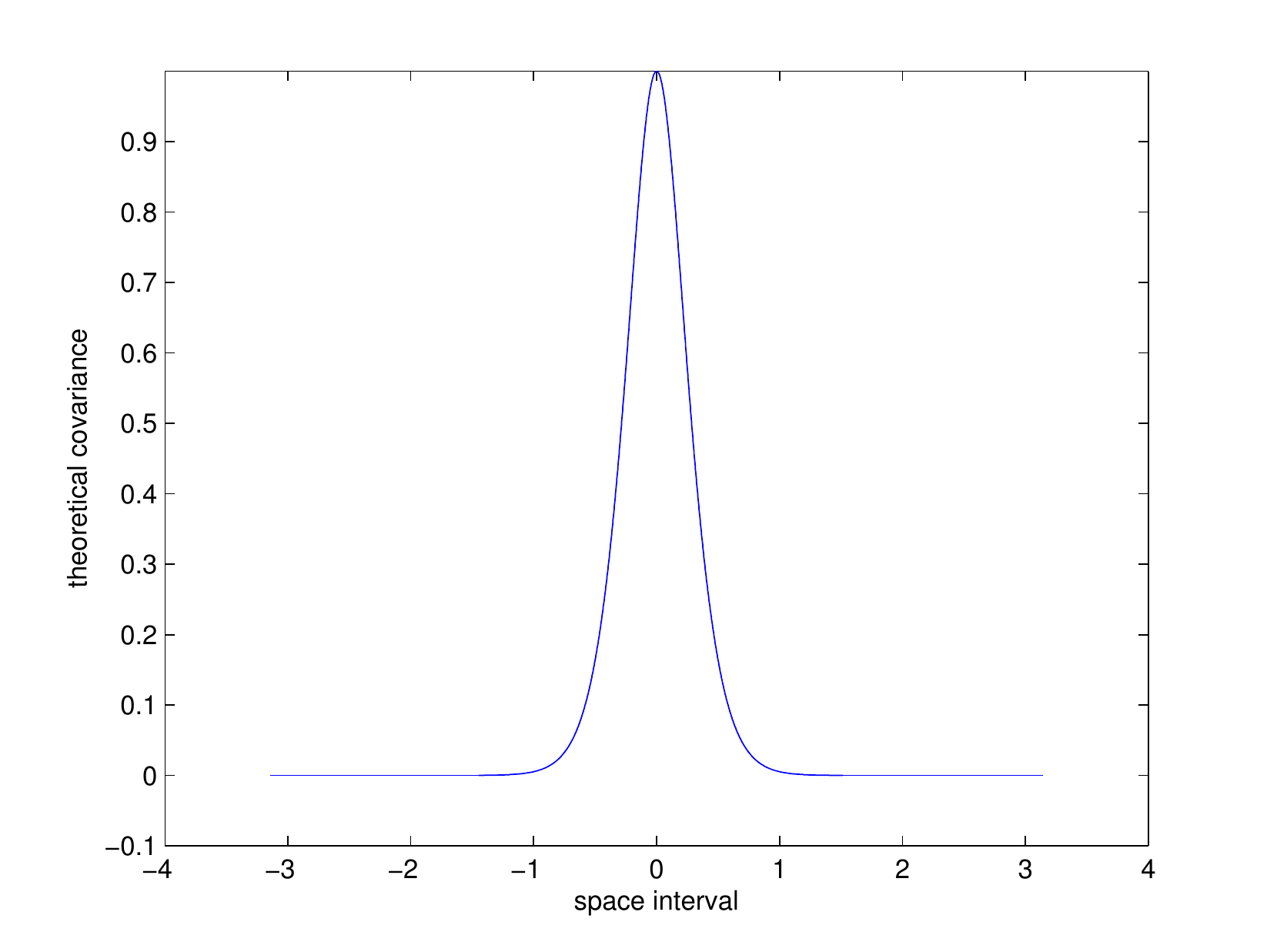}}
\subfigure[9 grid points]{
\includegraphics[width=0.48\textwidth]{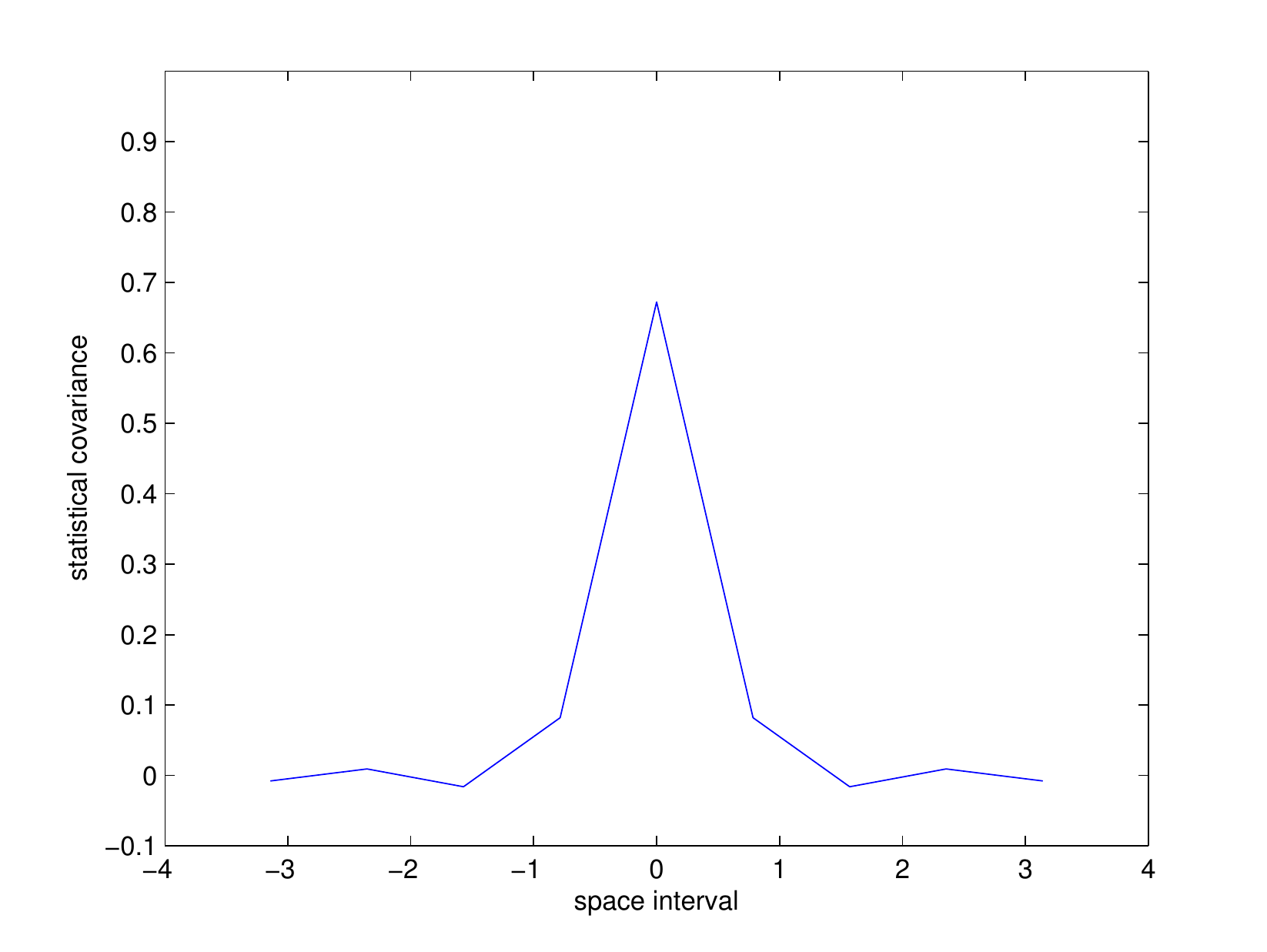}}
\subfigure[33 grid points]{
\includegraphics[width=0.48\textwidth]{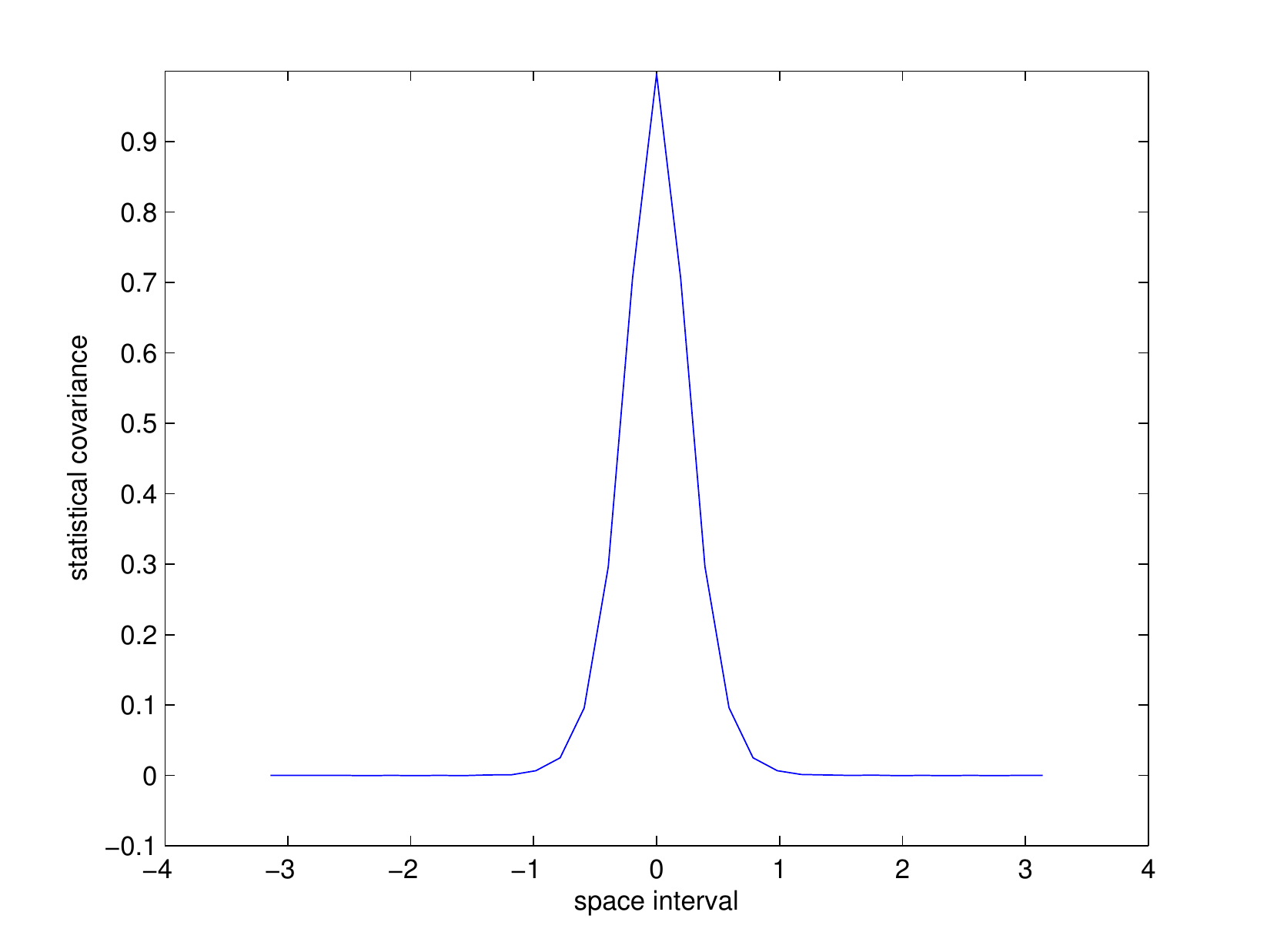}}
\subfigure[4097 grid points]{
\includegraphics[width=0.48\textwidth]{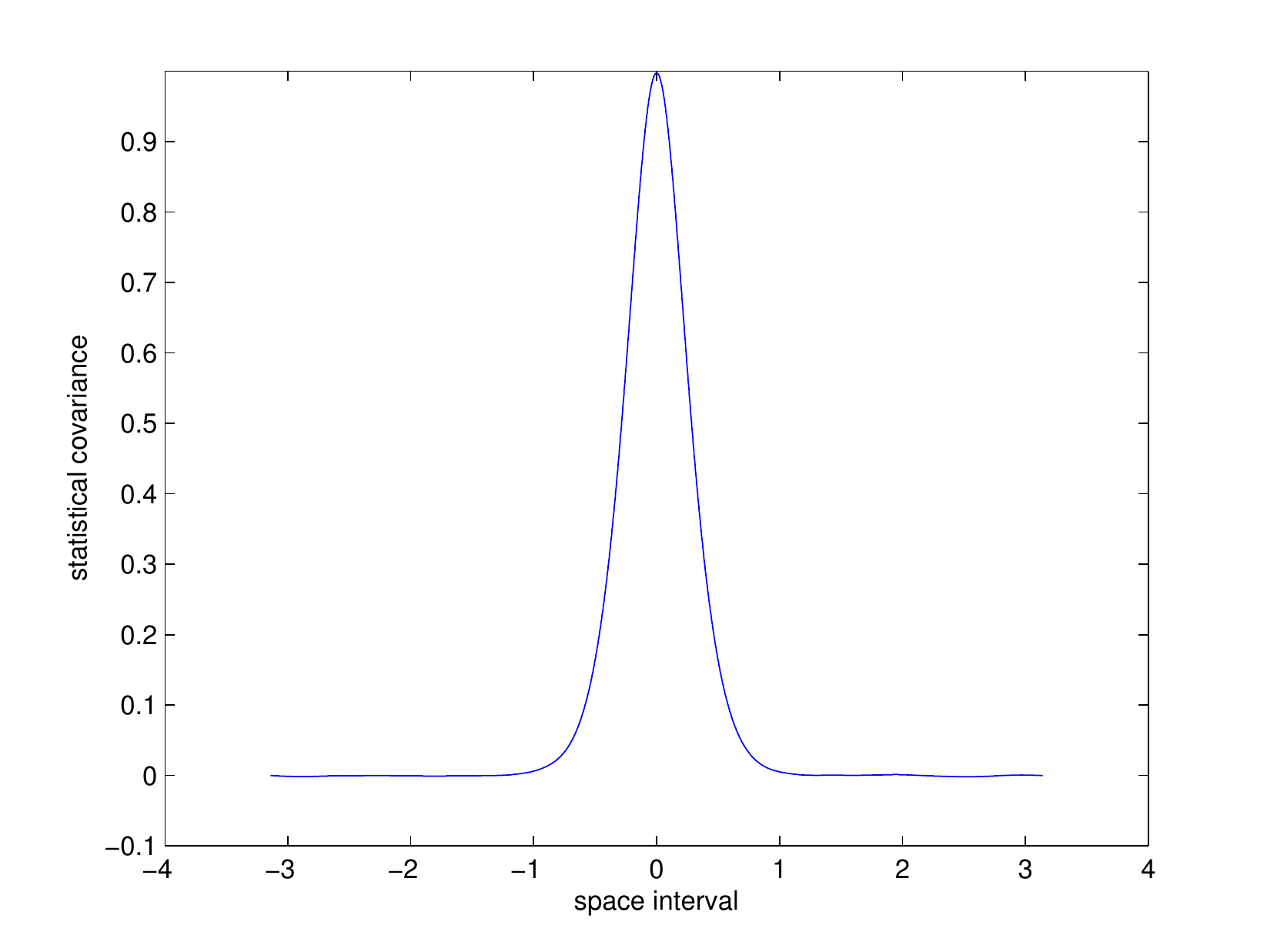}}
\caption{Theoretical and statistical covariance on
different grids.\label{figure_stat_cov_1d}}
\end{figure}

The error was computed by taking on each grid the maximum over all grid points of the
difference of the theoretical and statistical covariance, i.e. the error~$e^n$
on~$D^n$ was computed by
\begin{align*}
 e^n
  = \max_{i = -2^{n-1}, \ldots, 2^{n-1}}
      |\overline{\text{Cov}}(x_i^n) - C(x_i^n)|.
\end{align*}
The results for $n=2,\ldots,6$ are shown in Figure~\ref{figure_error_cov_1d}.
\begin{figure}
 \centering
\includegraphics[width=\textwidth]{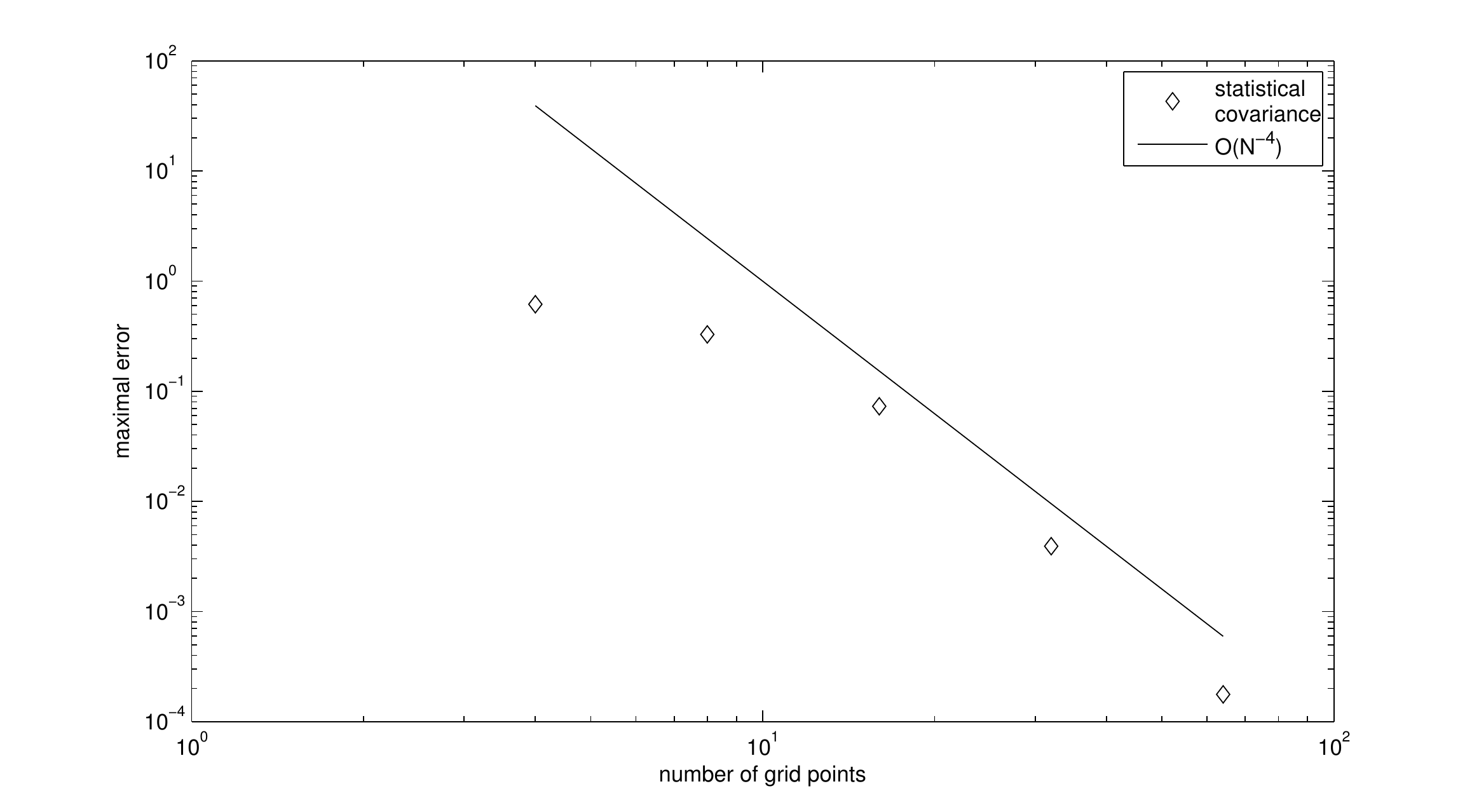}
\caption{Maximal difference of the statistical and theoretical covariance over all
grid points and the reference slope of $\Op(N^{-4})$.\label{figure_error_cov_1d}}
\end{figure}
The reference slope shows that the order of convergence is at least of
order~$\Op(N^{-4})$ where $N = |D^n|$. Finer grids are excluded from the plot
because then the error is dominated by the error of the Monte Carlo sampling.

\begin{appendix}
\section{Construction of $\cL^N$}
We suppose that $d\in\N$, and that $N=(N_1,\dotsc, N_d)$ with even $N_1$, \dots, $N_d\in\N$.
Put $M_i=N_i/2$, $i=1$, \dots,~$d$. Recall that
\begin{equation*}
	\cK^N = \bigl\{(k_1,\dotsc,k_d),\,k_i = 0,\dotsc, N_i-1,\,i=1,\dotsc, d\bigr\},
\end{equation*}
is subject to addition $\text{mod}(N)$. Set
\begin{equation*}
	\cL^N_0 = \bigl\{k\in\cK^N,\,k_i= \text{$0$ or $M_i$ for all
			    $i=1,\dotsc,d$}\bigr\}.
\end{equation*}
For $n=1$, \dots, $d$ construct $\cL^N_n$ as follows: Take all partitions
$(j_1,\dotsc,j_n)$, $(l_{n+1},\dotsc,l_d)$ of $(1,\dotsc,d)$, ordered in the
natural way, and collect $k\in \cK^N$ such that
\begin{align*}
	k_{j_1}	&= 1,\dotsc, M_{j_1}-1,\\
	k_{j_i}  &= 1,\dotsc, M_{j_i}-1,\,M_{j_i}+1,\dotsc, N_{j_i}-1,\quad i=1,\dotsc,n,\\
	k_{l_i}  &= 0,\,M_{l_i},\quad i=n+1,\dotsc, d.
\end{align*}
Then set
\begin{equation*}
	\cL^N = \biguplus_{n=0}^d \cL^N_n.
\end{equation*}
Observe that in the enumeration above, each $k\in\cL^N$ appears precisely
once. Moreover, the subset $\cK^N_0$ of $\cK^N$ so that $-\cK^N_0 = \cK^N_0$ coincides
with the subset $\cL^N_0$ of $\cL^N$. Set
\begin{equation*}
	\hat\cL^N = \cL^N\setminus\cL^N_0 = \biguplus_{n=1}^d \cL^N_n.
\end{equation*}
It is not hard to see that if $k\in\hat\cL^N$, then $-k\in\cK^N$ but
$-k\notin\cL^N$. Moreover, a simple counting argument shows that
$\cL^N\uplus (-\hat\cL^N)$ has the same number of elements as $\cK^N$.
Thus we find that $\cK^N = \cL^N_0\uplus\hat\cL^N \uplus(-\hat\cL^N)$,
and hence $\cL^N$ has indeed all required properties.

\end{appendix}

\providecommand{\bysame}{\leavevmode\hbox to3em{\hrulefill}\thinspace}
\providecommand{\MR}{\relax\ifhmode\unskip\space\fi MR }
\providecommand{\MRhref}[2]{%
  \href{http://www.ams.org/mathscinet-getitem?mr=#1}{#2}
}
\providecommand{\href}[2]{#2}

\end{document}